\documentclass{amsart}

\usepackage{a4}
\usepackage{latexsym}
\usepackage{amsmath}

\title[Generalised Difference Families]{All Difference Family Structures\\ arise from Groups}

\author{Tim Boykett}

\newtheorem{prop}{Proposition}

\newtheorem{lemma}[prop]{Lemma}
\newtheorem{thm}[prop]{Theorem}
\newtheorem{defn}[prop]{Definition}
\newtheorem{cor}[prop]{Corollary}

\newcommand{\beq}{\begin{equation}}
\newcommand{\eeq}{\end{equation}}

\newcommand{\dev}{dev}
\newcommand{\cal}{\mathcal}

\makeindex

\begin{document}

\bibliographystyle{plain}

\begin{abstract}
Difference families are traditionally built using groups as their
basis. This paper looks at what sort of 
generalised difference family
constructions could be made, using the standard basis
of translation and difference.

The main result is that  minimal requirements
on the  structure force  nothing that  groups cannot give,
at least in the finite case. Thus all
difference families arise from groups.
\end{abstract}

\maketitle

\section{Introduction}

Difference families are used to construct 2--designs.
They are based upon groups, usually additively written.
The essential operations are the difference operation
and the translation. 
The differences need to remain invariant under translation.

Quasigroups and loops are generalisations of groups that
do not require the operation to be associative.
A \emph{quasigroup} is a 2--algebra $(S,+)$ such that
for all $a,b \in S$ all equations
\[
a+x=b \;\;\; y+a=b
\]
have unique solutions for $x$ and $y$. The Cayley tables of such algebras
form Latin squares. There are many special cases of such
algebras. In particular a \emph{loop} is a quasigroup with
a two--sided identity $e \in S$.
A group is an associative loop. See for instance \cite{pflugfelder90}.

Quasigroups can be obtained by twisting a group
in some way. A simple example is to take an
additive group $(G,+)$ and to use the subtraction operation
 to obtain a quasigroup $(G,-)$ that is not associative.
For another example, given a field $K$, let $k\in K$, $k\neq 0,1$ be
arbitrary but fixed. Define
$a*b = ka+(1-k)b$. Then $(K,*)$ is a quasigroup, in general nonassociative.

There exists a more general form of equivalence between
quasigroups, or more general algebras. Two groupoids
$(S,+)$ and $(T,*)$ are \emph{isotopic} if there exist
bijections $\alpha,\beta,\gamma : S \rightarrow T$
such that for all $a,b\in S$
$\alpha (a+b) = \beta (a) * \gamma(b)$.
An isomorphism is an isotropism with all bijections
identical.

The isotope of a quasigroup is a quasigroup. In particular,
many quasigroups are isotopic to groups.

In the following we will first look at difference families and
determine what properties are needed to have in order
to be useful for such a construction. We demonstrate that
such a structure is equivalent to a class of
quasigroups.
We will then look at this
class of
algebras and see that they are all
simply obtained from groups.
The difference family
structures come directly from that group.
Our main results are Proposition \ref{prop_general} which
gives an explicit construction of all such quasigroups
and Proposition \ref{prop_ident} demonstrating that
the difference families are identical.

In general we are only interested in finite structures. However,
almost all the results here also apply for infinite structures.
We will note the use of finiteness arguments, 
which are only used in section \ref{sec_algprop}.

\section{Difference Families}

A \emph{(set) 2--design} is a pair $(V,{\cal B})$, where ${\cal B}$
is a set of subsets of $V$ all of size $k$ and for all pairs $a,b \in V,\,
a\neq b$, $|\{B \in {\cal B}:a,b \in B\}| = \lambda$ for 
some constants $k$ and $\lambda$.
The number 2 in the name refers to the pairs of elements $a,b$.
There are many variations on this definition, 
see e.g.\ \cite{Beth:DesignI} for details.

Given a (2)--algebra $(N,+)$ and a set of subsets $\cal B$ of $N$,
define \emph{development} of $\cal B$ in $N$ $dev({\cal B})$ 
to be the collection  $\{B+n: B \in {\cal B},\,n\in N\}$, possibly
containing duplicates. The set development is the collection
with no duplicates.

 Given a group $(N,+)$, not neccesarily abelian,
and a set ${\cal B} = \{B_i| i=1,\ldots,s \}$ of subsets of $N$, 
called \emph{base blocks}, such that
\begin{itemize}
\item all $B_i$ have the same size
\item for all $B,C\in {\cal B}$, 
$n\in N$, $B + n = C \Leftrightarrow B=C$ and $n=0$
\item  there exists some $\lambda$ such that for all nonzero $d \in N$,
$|\{(B,a,b) | B \in {\cal B},\,a,b \in B, a-b=d\}| = \lambda$
\end{itemize}
Then $N$ and $B$ form a \emph{difference family}, (DF). 

\begin{thm}[see e.g.\ \cite{Beth:DesignI}]
Let ${\cal B}$ be a difference family on a group $(N,+)$.
Then $dev {\cal B}$ is a set 2-design.
\end{thm}

In the proof of this
result, we can see that the requirement that $(N,+)$ be a group
is too strong. We are only using the translation property
 and the difference operation.
Thus it would seem that this construction can be generalised
to be based upon other structures. The following result does this.

\begin{thm}
\label{thm_very}
Let $N$ be a set with binary operation $-$ (difference) and unary operations
$t_i \in T$ (translations) and ${\cal B}$ a set of subsets of $N$ such that
\begin{itemize}
\item for all $a,b\in N$ there is a unique $t_i$ such that $t_i a = b$.
\item for all $a,b\in N$, the equation $a-x=b$ has a unique solution.
\item $a-b = t_ia - t_i b$ for all $a,b\in N$, for all $t_i$.
\item there exists some $\lambda$ such that for all  $d \in N$
such that $d=a-b$ for some $a,b \in N$,
$\Delta(d) = \{(B,a,b) |B \in {\cal B},\, a,b \in B,\, a\neq b,\, a-b=d\}$,
$|\Delta(d)| = \lambda$
\item there exists some integer $k$ such that $|t_i B| = k$ for all $t_i$ for
all $B \in {\cal B}$.
\item $t_i B = t_j C$ for $B,C \in {\cal B}$ implies $i=j$ and $B=C$.
\end{itemize}
Then $dev {\cal B} = \{t_i B: B \in {\cal B}\}$ is a set 2-design.
\end{thm}
\begin{proof}
All the blocks in $dev {\cal B}$ have size $k$ by construction.
They are all distinct by the final requirement. 
We need only show that the number of blocks
on a pair of points is constant.

Let $a\neq b \in N$. 
We show that 
$|\Delta(a-b)| = |\{t_i B: t_i\in T,B\in {\cal B},  a,b, \in t_iB\}|$.
There are exactly 
$\lambda$ triples $(B,\alpha,\beta)$ in $\Delta(a-b)$
such that $\alpha,\beta \in B$, $\alpha-\beta = a-b$.
For each $(B,\alpha,\beta)$ in $\Delta(a-b)$ 
there is a unique $t_i$ such that $t_i \alpha = a$.
We know
\[
a-b = \alpha - \beta = t_i \alpha - t_i \beta = a-t_i \beta
\]
so by the unique solution property of difference, $b=t_i \beta$.
Thus $a,b \in t_i B$, so we have a mapping $\Theta$ from
$\Delta(a-b)$ into $\{t_i B: t_i\in T,B\in {\cal B},  a,b, \in t_iB\}$. 
This map $\Theta$ is injective by the final condition.

We now show that $\Theta$ is surjective.
Let  $a,b \in t_i B$. Then
there exist some $\alpha, \beta \in B$ such that
$a = t_i \alpha$, $b=t_i \beta$, 
\[\alpha - \beta =  t_i \alpha - t_i \beta = a-b
\]
so $(B,\alpha,\beta) \in \Delta(a-b)$, and $t_i B$ is in
the image of $\Theta$. 
Thus $\Theta$ is a bijection and we are done.
\end{proof}

We now have a generalised form of difference family.
In the next sections we will investigate the algebraic properties
underlying this result.

\section{Algebraic Properties}
\label{sec_algprop}

Let us investigate the algebraic properties of the
results above. For this section, let $N$, $t_i$ and
$-$ be as defined in Theorem \ref{thm_very} above.

\begin{lemma}
\label{lemma_ident}
There exists a unique $t_0$ that is the identity
mapping on $N$.
\end{lemma}
\begin{proof}
Fix $a\in N$.
By the unique solution property of translations, there is
some $t_0\in T$ such that $t_0a = a$.
Then for all $b\in N$, 
$a-b = t_0a - t_0 b = a-t_0 b$, so by the unique
solution property of differences, $t_0 b = b$ and
$t_0$ is the identity map on $N$.
\end{proof}

\begin{thm}
If $N$ is finite then $(N,-)$ is a quasigroup.
\end{thm}
\begin{proof}
We know $a-x=b$ has a unique solution. Suppose
$x-a=b$ has  two solutions, $x_1\neq x_2$.
There is some $t_i$ such that $t_ix_1 = x_2$.
Then 
\[
x_1 - a = t_ix_1 - t_ia = x_2 - t_ia = x_2 - a = b
\]
Then $x_2 - x = b$ has solutions $a$ and $t_ia$ for $x$.
Thus $t_ia = a$ so $t_i$ is the identity and $x_1 = x_2$ thus
$x-a = b$ has at most one solution. By finiteness, it has
exactly one solution, so $(N,-)$ is a quasigroup.
\end{proof}

\begin{defn}
Fix $a_0\in N$. For all $b\in N$, let $t_b \in T$ such that
$t_b a_0 = b$. This is unique. Define 
$x+y := t_y x$.
\end{defn}

\begin{thm}
If $N$ is finite then $(N,+)$ is a quasigroup.
\end{thm}
\begin{proof}
By the uniqueness property above, $a+x=b$ has a unique solution.
Suppose $x+a=b$ has two solutions $x_1, x_2$.
Then
$x_1 - x_2 = (x_1+a) - (x_2+a) = b-b$.
There is some unique $k$ such that $b+k=x_1$.
Then 
$x_1 - x_2 = b-b = (b+k) - (b+k) = x_1 - x_1$ so
$x_1 = x_2$ by quasigroup property of $(N,-)$,
so $x+a=b$ has at most one solution. By finteness
it has exactly one, and $(N,+)$ is a quasigroup.
\end{proof}

Thus we have shown that the structure used in Theorem
\ref{thm_very} can be seen as a set with two
operations that form quasigroups. We formalise this,
as it is clear that from such a pair of quasigroups,
we can form the translations used in Theorem \ref{thm_very}.

\begin{defn}
A \emph{difference family biquasigroup (DFBQ)} $(N,+,-)$ is a $(2,2)$--algebra
where each operation gives a quasigroup and the
equation $a-b = (a+c)-(b+c)$ is satisfied.
\end{defn}

\begin{lemma}
A DFBQ has a right additive identity.
There is a constant $e\in N$ such that  $e=a-a$ for all $a$.
\end{lemma}
\begin{proof}
Let $a,\bar a \in N$ such that $a+\bar a = a$.
Then for any $b$, 
\[
b-a = (b+\bar a) - (a+\bar a) = (b+\bar a) - a.
\]
By the quasigroup property, $b=b+\bar a$.
This works for all $a,b$, so $\bar a$ is constant,
and this constant is a right identity with respect to
addition. This identity is unique by the quasigroup
property.

Fix some $a\in N$. Define $e:=a-a$. For all $b\in N$, there
exists some $c$ such that $a+c=b$. Thus
$b-b = (a+c)-(a+c) = a-a = e$ and the second statement
is proved.
\end{proof}

Note that the first part also follows from Lemma \ref{lemma_ident} above.
We may write a DFBQ as $(N,+,-,o,e)$ where
$o$ is a right additive identity and $e=a-a$ for all $a$.

\section{In general}

In this section we examine the structure of
a general DFBQ. We will use these results in
the next section to demonstrate that a general
DFBQ is isotopic to a group and that
the resulting designs are identical.

\begin{defn} A collection ${\cal B}$ and a DFBQ $(N,+,-)$
such that:
\begin{itemize}
\item there exists an integer $k$ such that $|B| = k$ for all $B \in {\cal B}$
\item there exists some $\lambda$ such that for all  $d \in N$
such that $d=a-b$ for some $a,b \in N$,
$\Delta(d) = \{(B,a,b) |B \in {\cal B},\, a,b \in B,\, a\neq b,\, a-b=d\}$,
$|\Delta(d)| = \lambda$
\item $ B+b = C+c$ for $B,C \in {\cal B}$, $b,c\in N$ implies $B=C$ and $b=c$.
\end{itemize}
 is called a \emph{quasigroup difference family (QDF)}
\end{defn}

It is clear that such a QDF will give a 2-design using the same
methods as Theorem \ref{thm_very}.

\begin{prop}
\label{prop_breakdown}
Let $(N,+,-,o,e)$ be a DFBQ. Let $\bar e$ be such that $e+\bar e = o$,
then define $\phi:x\mapsto x+\bar e$ and $\alpha: x \mapsto x-e$.
Define the operations
\begin{eqnarray}
a \oplus b &=& \phi^{-1}(\phi a + \phi b) \\
a \ominus b &=& \alpha^{-1}(a-b)
\end{eqnarray}
Then $(N,\oplus,\ominus,e,e)$ is a DFBQ with $a\ominus e = a$ for all
$a\in N$.
\end{prop}
\begin{proof}
$(N,\oplus)$ and $(N,\ominus)$ are quasigroups by isotopism.
Note that $\phi a - \phi b = a-b$ by the DFBQ property,
so $\phi^{-1} a - \phi^{-1} b = a-b$.
Then
\begin{eqnarray}
(a\oplus c) \ominus (b\oplus c) &=& \alpha^{-1}((a\oplus c) - (b\oplus c)\\
 &=& \alpha^{-1}((\phi a + \phi c)-(\phi b + \phi c))\\
  &=& \alpha^{-1}(\phi a - \phi b) = \alpha^{-1}(a-b)\\
 &=& a \ominus b
\end{eqnarray}
so we have the DFBQ property.
The constants are both $e$ as given by
$a\oplus e = \phi^{-1}(\phi a + \phi  e) = \phi^{-1}(\phi a + o) = a$
and $a\ominus a =  \alpha^{-1}(a-a) = \alpha^{-1}(e)=e$.

The second claim is seen by $\alpha a = a-e$ thus
$a\ominus e = \alpha^{-1}(a-e) = a$.
\end{proof}
Note that $a\ominus b = \phi a \ominus \phi b$. This will be important for the
next result.

\begin{prop}
\label{prop_backup}
Let $(N,\oplus,\ominus,e,e)$ be a DFBQ with $a\ominus e = a$ for all
$a\in N$, $\phi$ a permutation of $N$ 
such that $a\ominus b = \phi a \ominus \phi b$, $\alpha$ a
permutation of $N$ such that $\alpha e = e$.
Define
\begin{eqnarray}
a +b &=& \phi(\phi^{-1} a \oplus \phi^{-1} b)\\
a-b &=& \alpha (\phi^{-1} a \ominus \phi^{-1} b)
\end{eqnarray}
Let $o := \phi e$.
Then $(N,+,-,o,e)$ is a DFBQ,
$\alpha:a\mapsto a-e$, $\phi: a \mapsto a+\bar e$ where $e+\bar e = o$.
\end{prop}
\begin{proof}
$(N,+)$ and $(N,-)$ are quasigroups by isotopism.
The DFBQ property is seen by
\begin{eqnarray}
(a+c)-(b+c) &=& 
\alpha(\phi^{-1}\phi(\phi^{-1}a \oplus \phi^{-1}c)\ominus
       \phi^{-1}\phi(\phi^{-1}b \oplus \phi^{-1}c))\\
&=& \alpha((\phi^{-1}a \oplus \phi^{-1}c)\ominus
           (\phi^{-1}b \oplus \phi^{-1}c))\\
&=& \alpha(\phi^{-1}a \ominus    \phi^{-1}b )\\
&=& a-b
\end{eqnarray}
The constants are given by
$a+o = \phi(\phi^{-1} a \oplus \phi^{-1}o) =  \phi(\phi^{-1} a \oplus e) = a$
and $a-a = \alpha (a\ominus a) = \alpha e = e$.

Then
$\alpha(a) = \alpha(a\ominus e) = \alpha (\phi^{-1}a \ominus \phi^{-1}e)
= a-e$.
Let $\bar e$ be such that $e+\bar e = o$.
Then $e+\bar e = \phi(\phi^{-1}e \oplus \phi^{-1}\bar e) = o = \phi e$
so $\phi^{-1}e \oplus \phi^{-1}\bar e = e$.
Then
\begin{eqnarray}
a=a\ominus e &=&  \phi^{-1} a \ominus\phi^{-1}e\\
&=& ( \phi^{-1} a \oplus \phi^{-1}\bar e) \ominus (\phi^{-1}e \oplus \phi^{-1}\bar e)\\
&=& ( \phi^{-1} a \oplus \phi^{-1}\bar e) \ominus e \\
&=& \phi^{-1} a \oplus \phi^{-1}\bar e\\
&=& \phi^{-1}(a+\bar e)
\end{eqnarray}
Thus $\phi a = a+\bar e$ and we are done.
\end{proof}

We will need the following for the final results.

\begin{defn}
A quasigroup $(Q,\circ)$ is a \emph{Ward quasigroup}
if $(a\circ c)\circ (b\circ c) = a\circ b$ for all $a,b,c\in Q$.
\end{defn}

\begin{thm}[\cite{vojtechovsky03}]
Let $(Q,\circ)$ be a Ward quasigroup. Then there exists a unique
element $e\in Q$ such that for all $x\in Q$, $x\circ x=e$.
Define $\bar x = e\circ x$ and $x* y = x\circ\bar y$ for all $x,y \in Q$.
Then $(Q,*,\bar{})$ is a group, and $x\circ y=x* \bar y$.
\end{thm}

\begin{prop}
\label{propuseward}
Let $(N,+,-,e,e)$ be a DF biquasigroup with $a-e=a$ for all $a$.
Then it is isotopic to a DF group.
\end{prop}
\begin{proof}
Let $I$ be the permutation of $N$ such that $a+Ia=e$.
Then 
\begin{equation}
a-b = (a+Ib)-(b+Ib) = (a+Ib) - e = a+Ib. \label{eq_useward}
\end{equation}
Thus $(a-b)-(c-b) = (a+Ib)-(c+Ib) = a-c$ so $(N,-)$ is
a Ward quasigroup.
Thus there is a group $(N,*,\cdot^{-1})$ with 
$a-b = a*b^{-1}$ and $a+b=a*(I^{-1}b)^{-1}$ 
by equation (\ref{eq_useward}).
\end{proof}

The converse of this result holds too. The proof is simple calculation.

\begin{lemma}
\label{lemma_Ifix1}
Let $(N,*,1)$ be a group, $I$ a permutation of $N$ fixing $1$.
Define
\begin{eqnarray}
a+b &=& a*(Ib)^{-1}\\
a-b &=& a*b^{-1}
\end{eqnarray}
Then $(N,+,-,1,1)$ is a DF biquasigroup with $a-1=a$ for all $a$.
\end{lemma}

Thus we obtain information on the form of the map $\phi$ in Proposition
\ref{prop_backup}. We know the form of the operations from Prop
\ref{propuseward} so we can make some explicit statements about
the structure.

\begin{cor}
\label{cor_phi}
Let $(N,\oplus,\ominus,e,e)$ and $\phi$ 
be as for Proposition \ref{prop_backup}.
Let the operation $*$ be as from Prop \ref{propuseward}.
Then there exists some $k\in N$ such that the 
map $\phi$  is of the form
$\phi(a) = a*k$.

Conversely, given  $(N,\oplus,\ominus,e,e)$ as in Prop \ref{prop_backup}
and a group operation $*$, select any element $k\in N$.
Then $\phi(a) := a*k$ satisfies the requirements of Prop \ref{prop_backup}.
\end{cor}
\begin{proof}
By Prop \ref{propuseward} we know that 
$a\ominus b = a*b^{-1}$. Since $\phi a \ominus \phi b = a \ominus b$
we have
$\phi a * (\phi b)^{-1} = a * b^{-1}$. Let $b=1$ and we
obtain
$\phi a *(\phi 1)^{-1} = a$ so $\phi a = a*\phi 1$. Letting $k:= \phi 1$
we are done.

The converse is seen by taking any element $k\in N$.
Define $\phi a := a*k$.
Then $\phi a \ominus \phi b = (a*k) * (b*k)^{-1} = a*b^{-1} = a \ominus b$
so we are done.
\end{proof}

\section{General Explicit Descriptions}

In this section, we will look at explicit descriptions of
DFBQs and QDFs.
Using the results above,
we know the structure of all DFBQs.

\begin{prop}
\label{prop_general}
Let $(N,*,1)$ be a group. Let $\alpha,\beta$ be permutations
of $N$,  $\alpha 1 = 1$. Define
\begin{eqnarray}
a+b &=& a * \beta b\\
o &=& \beta^{-1}(1)\\
a-b &=& \alpha(a* b^{-1})
\end{eqnarray}
Then  $(N,+,-,o,1)$ is a DFBQ and all  DFBQs are 
of this form.
\end{prop}

\begin{proof}
The forward direction is a calculation and is clear.
Let $(N,+,-,o,e)$ be a DFBQ. We demonstrate that
there exists a group structure $(N,*,{}^{-1},1)$ and
 permutations $\alpha,\beta$ of $N$ as above.

By Proposition \ref{prop_breakdown}
there exist $\phi$ and $\alpha$ such that defining
\begin{eqnarray}
a \oplus b := \phi^{-1} (\phi a + \phi b)\\
a \ominus b := \alpha^{-1} (a-b)
\end{eqnarray}
we obtain $(N,\oplus,\ominus,e,e)$ is a DFBQ with
$a\ominus e = e$.
By Proposition \ref{propuseward} there exists some group
$(N,*,{}^{-1},1)$
such that $e=1$, $a\ominus b = a*b^{-1}$ and $a\oplus b = a*(I^{-1}(b))^{-1}$.
Thus
\begin{eqnarray}
a+b &=& \phi(\phi^{-1}a * (I^{-1}(\phi^{-1}b))^{-1})\\
a-b &=& \alpha(a*b^{-1})
\end{eqnarray}

By Corollary \ref{cor_phi} we know that $\phi x = x* k$, 
$\phi^{-1}x = x*k^{-1}$.
Thus
\begin{eqnarray}
a+b &=& ((a*k^{-1})*(I^{-1}(b*k^{-1}))^{-1})*k \\
  &=& a * k^{-1}*(I^{-1}(b*k^{-1}))^{-1}*k \\
  &=& a* \beta (b)
\end{eqnarray}
where $\beta (x) = k^{-1}*(I^{-1}(x*k^{-1}))^{-1}*k$ is a permutation
of $N$.
Since $a+\beta^{-1}(1) = a*\beta(\beta^{-1}(1)) = a*1=a$ we know
$\beta^{-1}(1)$ is the unique right identity, so 
$o=\beta^{-1}(1)$. The permutation $\alpha$ fixes $e$ which
is seen to be $1$ and we are done.
\end{proof}

This final result shows that all difference family structures
are in fact group structures.

\begin{prop}
\label{prop_ident}
The quasigroup development and the group
development of a difference family are identical.
\end{prop}
\begin{proof}
Suppose we have a QDF ${\cal B}$ on a DFBQ $(N,+,-,o,e)$.
By Prop \ref{prop_general} above, we know that there
is a group operation $*$ and some permutation of $N$ such that
$a+b = a*\beta(b)$. Thus if $B$ is a subset of
$N$, 
\[
dev_+ B = \{B+n: n\in N\} = \{B*\beta(n) : n\in N\} =  \{B*n : n\in N\} = dev_* B
\]
so we obtain exactly the same set of sets.
Thus $dev_+ {\cal B} = \dev_*{\cal B}$ and we are done.
\end{proof}

\section{Conclusion}
It would be desirable to generalise the definition of a
difference family so as to use more general  structures
to derive designs using this formalism.
With simple and reasonable requirements for our difference family
structures, we have shown that we obtain a
biquasigroup algebra and that such algebraic
structures must be isotopic to groups. It is also seen
that the resulting designs are identical.

Questions remain open as to whether the requirements that we posit
are all necessary. It may be reasonable to use a simpler
structure for the difference operation, but I cannot see how.

Applications remain open here. For instance, planar nearrings
have been shown to possess a difference family structure.
Questions about nonassociative planar nearrings have been 
raised, and it might be appropriate to use these results
to deduce structure about the nearrings that could
be so defined. It also remains open as to the 
properties of infinite generalised difference families, where
the translations and difference operation do not form a proper quasigroup.
The investigation of neardomains and K--loops \cite{kiechle02} suggests
that there are some strange and interesting properties when we drop the
finiteness restriction. In particular there may be 
connections between the generalisation of nearfields to neardomains
and the generalisation to planar nearrings and Ferrero pairs
\cite{clay92,pilz83},
which may be connected to the construction of
nonassociative difference families.

\section{Acknowledgements}

This work has been supported by grant P15691 of 
the Austrian National Science Foundation 
(Fonds zur F\"orderung der wissenschaftlichen Forschung).
I would also like to thank Petr Vojtechovsky for many
ideas and support with loops and quasigroup
theory.

\bibliography{tims}

\end{document}